\documentclass[final]{siamltex}
\pdfoutput=1 

\usepackage{epsfig}
\usepackage{times}
\usepackage{rotating}
\usepackage{xspace}
\usepackage{graphicx}

 \def\ul{\underline}

\def\eref#1{(\ref{#1})}

\def\Proof{\noindent {\bf Proof.  }}

\newtheorem{THM}{THEOREM}[section] \newtheorem{LEMMA}[THM]{Lemma}
\newtheorem{CLAIM}[THM]{Claim} \newtheorem{COR}[THM]{Corollary}
\newtheorem{PROP}[THM]{Proposition}

\newcommand{\eq}{\begin{equation}} \newcommand{\eeq}{\end{equation}}

\newcommand{\THEOREM}{\begin{THM}} \newcommand{\eT}{\end{THM}}
\newcommand{\Lemma}{\begin{LEMMA}} \newcommand{\eL}{\end{LEMMA}}
\newcommand{\Claim}{\begin{CLAIM}} \newcommand{\eCl}{\end{CLAIM}}
\newcommand{\Corollary}{\begin{COR}} \newcommand{\eCo}{\end{COR}}

\newcommand{\Proposition}{\begin{PROP}} \newcommand{\eP}{\end{PROP}}

\title{Sink-Stable Sets of Digraphs}

\author{D\'ora Erd{\H o}s\thanks{ Department of Computer Science,
Boston University 111 Cummington St, Boston MA 02215, e-mail:  {\tt
edori\char'100 cs.bu.edu .} } \and Andr\'as Frank,\thanks{MTA-ELTE
Egerv\'ary Research Group, Department of Operations Research,
E\"otv\"os University, P\'azm\'any P. s. 1/c, Budapest, Hungary,
H-1117. e-mail:  {\tt frank\char'100 cs.elte.hu .}} \and Kriszti\'an
Kun\thanks{Webra International Kft, Andocs u. 45, Budapest, Hungary,
H-1165 e-mail:  {\tt kknc\char'100kethold.hu .} } }

\begin{document}

\maketitle

\begin{abstract} We introduce the notion of sink-stable sets of a
digraph and prove a min-max formula for the maximum cardinality of the
union of $k$ sink-stable sets.  The results imply a recent min-max
theorem of Abeledo and Atkinson \cite{Abeledo-Atkinson1} on the Clar
number of bipartite plane graphs and a sharpening of Minty's coloring
theorem \cite{Minty62}.  We also exhibit a link to min-max results of
Bessy and Thomass\'e \cite{Bessy-Thomasse} and of Seb{\H o}
\cite{Sebo07} on cyclic stable sets.

\end{abstract}

\section{Introduction}

It is well-known that the problem of finding a stable set of maximum
cardinality is {\bf NP}-complete in a general undirected graph but
nicely tractable, for example, for comparability graphs.  A
comparability graph is the underlying undirected graph of a
comparability digraph, which is, by definition, an acyclic and
transitive digraph.  Such a digraph can also be considered as one
describing the relations between the pair of elements of a partially
ordered set.  A subset $S$ of nodes of a directed graph $D=(V,A)$ is
defined to be stable if $S$ is stable in the underlying undirected
graph of $D$.  Dilworth's theorem \cite{Dilworth50} states in these
terms that the maximum cardinality of a stable set of a comparability
digraph is equal to the minimum number of cliques covering $V$.
Greene and Kleitman \cite{Greene-Kleitman} extended this result to a
min-max theorem on the maximum cardinality of the union of $k$ stable
sets of a comparability digraph.

The present investigations have three apparently unrelated sources.
In solving a long-standig conjecture of Gallai \cite{Gallai63}, Bessy
and Thomass\'e \cite{Bessy-Thomasse} introduced a special type of
stable sets, called cyclic stable sets, and proved a min-max result on
the maximum cardinality of a cyclic stable set.  They also derived a
theorem on the minimum number of cyclic stable sets required to cover
all nodes.  These two results were unified and extended by Seb{\H o}
\cite{Sebo07} who proved a min-max formula for the the largest union
of $k$ cyclic stable sets.  His theorem is an extension of the theorem
of Greene and Kleitman \cite{Greene-Kleitman}.  Another source is a
recent min-max result of Abeledo and Atkinson on the Clar number of
plane bipartite graphs.  The third source is a colouring theorem of
Minty \cite{Minty62}.  It will be shown that these remote results
have, in fact, a root in common.

To this end, we introduce and study another special kind of stable
sets of an arbitrary digraph $D$.  A {\bf directed cut} or {\bf dicut}
of a digraph is a subset of edges entering a subset $Z$ of nodes
provided no edge leaves $Z$.

A node of $D$ will be called a {\bf sink node} (or just a sink) if it
admits no leaving edges.  A node is a {\bf source node} if it admits
no entering edge.  A subset of nodes is a {\bf sink set} if each of
its elements is a sink.  Clearly, a sink set is always stable.  We say
that a subset $S$ of nodes of $D$ is {\bf sink-stable} if there are
edge-disjoint directed cuts of $D$ so that reorienting the edges of
these dicuts $S$ becomes a sink set.  Obviously, any subset of a
sink-stable set is also sink-stable.  A source-stable set is defined
analogously.  Observe that a subset $S\subseteq V$ is sink-stable if
and only if $S$ is source-stable.  Note that a node of a directed
circuit never belongs to a sink-stable set since a dicut and a
di-circuit are always disjoint.  In the acyclic digraph with node-set
$\{a,b,c,d\}$ and edge-set $\{ab,\ bc,\ cd, \ ad\}$ every single node
is a one-element sink-stable set while the stable set $\{a,c\}$ is not
sink-stable.  Note that $D$ is not transitive and hence $D$ is not a
comparability digraph.

\begin{figure}[h]
\begin{center}
{\includegraphics[scale = 0.5,angle=270]{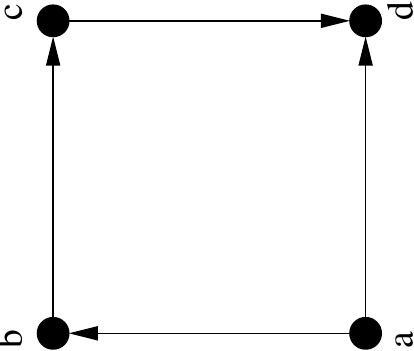}}%
\caption{\label{fig:ex_two_stable} Every node is a sink-stable set of size 1,
  while the stable set $\{a,c\}$ is not sink-stable.}
\end{center}
\vspace{-3ex}
\end{figure}

\Proposition \label{compar} In a comparability digraph $D=(V,A)$,
every stable set is sink-stable.  \eP

\Proof Let $S$ be a stable set.  We may assume that $S$ is maximal.
Let $Z$ denote the set of nodes of $V-S$ that can be reached along a
dipath from $S$.  We claim that no edge can leave $Z$.  Indeed, if
$uv$ does, then $v$ is also reachable form $s$ and hence $v$ must be
in $S$.  Since $D$ is acyclic, $v\not =s$.  Since $D$ is transitive,
there is an edge $sv\in A$, contradicting the stability of $S$.  It
follows that the edges entering $Z$ form a dicut and reorienting this
dicut $S$ becomes a sink set.  $\bullet $ \medskip

In the sequel, we shall also use a strongly related notion of special
stability.  For a subset $F$ of edges of a digraph $D$, a subset $S$
of nodes is {\bf $F$-stable} if $S$ is sink-stable in the digraph
$D_F$ arising from $D$ by reversing the elements of $F$.  It can be
checked that $S$ is sink-stable in $D$ if and only if $S$ is
$F$-stable in the digraph arising from $D$ by adding the reverse of
each edge of $D$ where $F$ denotes the set of these reversed edges.
It will turn out that in some cases it is easier to work out a result
for sink-stable sets and use then this to derive the corresponding
result for $F$-stable sets.  For example, in characterizing
sink-stable and $F$-stable sets we shall follow this path.  There are
other cases when the reverse approach is more convenient.  For
example, we shall prove first a min-max formula for $F$-stable sets
and use then this to derive the corresponding min-max theorem for
sink-stable sets.

The property of sink-stability is in {\bf NP} in the sense that the
set of disjoint dicuts whose reorientation turns a subset $S$ into a
sink set is a fast checkable certificate for $S$ to be sink-stable.
Theorem \ref{sinkjel} will describe a co-{\bf NP} characterization for
sink-stability.  We shall also characterize for any integer $k\geq 2$
the union of $k$ sink-stable sets, and as a main result, a min-max
formula will be proved for the largest union of $k$ sink-stable sets.

The result for $k=1$ shall imply a recent min-max theorem of Abeledo
and Atkinson \cite{Abeledo-Atkinson1} on the Clar number of a
2-connected bipartite plane graph $G$.  Here the Clar number is
defined to be the maximum number of disjoint bounded faces of $G$
whose removal leaves a perfectly matchable graph.  This notion was
originally introduced in chemistry for hexagonal plane graphs to
capture the behaviour of characteristic chemical and physical
properties of aromatic benzenoids.

We will also derive a sharpening of Minty's colouring theorem
\cite{Minty62} by proving a min-max formula for the minimum number of
sink-stable sets to cover $V$, and show how this result implies a
theorem of Bondy \cite{Bondy76} stating that the chromatic number of a
strongly connected digraph is at most the length of its largest
directed circuit.

Finally, an interesting link will be explored to a recent min-max
theorem of Bessy and Thomass\'e \cite{Bessy-Thomasse} on so-called
cyclic stable sets of strongly connected digraphs, a result that
implied a solution of a conjecture of Gallai.  A min-max theorem of
Seb{\H o} \cite{Sebo07} on the largest union of $k$ cyclic stable sets
will also be a consequence.


\medskip To conclude this introductory section, we introduce some
definitions and notation.  For a function $m:V\rightarrow {\bf R} $
(or vector $m\in {\bf R}\sp V$), we define a set-function $\widetilde
m$ by $ \widetilde m(X)= \sum [m(v):v\in X]$ where $X\subseteq V.$ For
a number $x$, let $x\sp +:=\max \{x, 0 \}$.  By a {\bf multi-set} $Z,$
we mean a collection of elements of $V$ where an element of $V$ may
occur in more than one copy.  The indicator function $\underline \chi
_Z:V\rightarrow \{0,1,2,\dots \}$ of $Z$ tells that $\underline \chi
_Z(v)$ copies of an element $v$ of $V$ occurs in $Z$.  A multiset is
sometimes identified with its indicator function which is a
non-negative integer-valued function on $V$.

Let $D=(V,A)$ be a digraph.  For function $x:A\rightarrow {\bf R} $,
the in-degree and out-degree functions $\varrho _x$ and $\delta _x$
are defined by $\varrho _x(Z)= \sum [x(uv):  uv\in A, \ u\in V-Z, \ v
\in Z]$ for $Z\subseteq V$ and by $\delta _x(Z):=\varrho _x(V-Z)$.\
$x$ is a {\bf circulation} if $\varrho _x=\delta _x$.  A function $\pi
:V\rightarrow {\bf R} $ is often called a potential.  For a potential
$\pi $, the {\bf potential difference} $\Delta _\pi :A\rightarrow {\bf
R} $ is defined by $\Delta _\pi (uv):=\pi (v)-\pi (u)$ where $uv\in
A$.  A function arising in this way is called a tension.

By a {\bf topological ordering} of a digraph, we mean an ordering
$v_1,\dots ,v_n$ of the nodes so that every edge $a$ of $D$ goes
forward, that is, $a$ is of type $v_iv_j$ where $i<j$.

A {\bf circuit} is a connected undirected graph in which the degree of
every node is 2. Typically, we use the convention for a circuit $C$
that $C$ also denotes the edge-set of the circuit while $V(C)$ denotes
its node-set.  A directed graph is also called a circuit if it arises
from an undirected circuit by arbitrarily orienting its edges.  Every
circuit $C$ of at least three nodes has two ways to traverse its
elements.  For a graph or digraph $H$ on node-set $V$ we fix an
ordering of the elements of $V$.  For every circuit $C$ of $H$ with at
least three nodes, the three smallest indexed nodes of $C$ uniquely
determine a traversal direction of $C$ called {\bf clockwise
direction} of $C$.  When $H$ is a digraph, the clockwise edges of $C$
will also be called {\bf forward} edges while the anti-clockwise edges
of $C$ are the {\bf backward} edges.

In a directed graph $D=(V,A)$, by a {\bf walk} $W$ we mean a sequence
$(v_0, e_1, v_1, e_2,\dots , e_k, v_k)$ consisting of not necessarily
distinct nodes and edges where $e_i$ is either a $v_{i-1}v_i$-edge
(called {\bf forward} edge) or a $v_iv_{i-1}$-edge (called {\bf
backward} edge).  If every edge is forward, we speak of a {\bf
one-way} walk.  If $v_0=v_k$, we speak of a {\bf closed walk}.  If the
terms of a closed walk are distinct apart form $v_0$ and $v_k$, we
speak of a simple closed walk.  Therefore a simple closed walk with at
least one edge can be identified with a circuit having a specified a
traversal direction.  Note that the closed walk consisting of a single
node and no edge is not a circuit.

The number of forward and backward edges of a circuit $C$ of a digraph
are denoted by $\varphi (C)$ and $\beta (C)$, respectively, while
their minimum will be called the {\bf $\eta $-value} of $C$ or
sometimes simply the {\bf value} of the circuit.  The value of $C$ is
denoted by $\eta (C).$ When $\eta (C)=0$, we speak of a {\bf
di-circuit}.  We emphasize the difference between a circuit whose
edges are just directed edges and a di-circuit.  An edge of a digraph
is {\bf cyclic} if it belongs to a di-circuit.  For a function
$x:A\rightarrow {\bf R}$, $\varphi _x(C)$ denotes the sum of the
$x$-values over the forward edges of circuit $C$ while $\beta _x(C)$
is the sum of the $x$-values over the backward edges.  Clearly,
$\varphi _x(C) + \beta _x(C)= \widetilde x(C)$.  For a subset $B$ of
edges, $\varphi _B(C)$ denotes the number of forward edges of $C$
belonging to $B$, while $\beta _B(C)$ is the number of backward edges
of $C$ belonging to $B$.

A function $x:A\rightarrow {\bf R} $ is {\bf conservative} if
$\widetilde c(K)\geq 0$ for every di-circuit $K$.  A potential $\pi $
is {\bf $c$-feasible} or just feasible if $\Delta _\pi \leq c$.

\medskip

\section{Dicut equivalence and sink-stable sets}

\bigskip

\Lemma [Gallai, \cite{Gallai58}] \label{konzerv} A cost function $c:A
\rightarrow {\bf R}$ on the edge-set of a digraph $D=(V,A)$ is
conservative if and only if there is a feasible potential.  Moreover,
if $c$ is integer-valued, then $\pi $ can also be selected to be
integer-valued.  $\bullet $ \eL

The lemma immediately implies for an integer-valued tension $x$ that
there is an integer-valued potential $\pi $ for which $x=\Delta _\pi
$.

\Lemma \label{dicuts} For a subset $F\subseteq A$ of edges of a
digraph $D=(V,A)$, the following are equivalent.  \medskip

{\em (A)} \ $F$ is the union of disjoint dicuts.  \medskip

{\em (B)} \ $\varphi _F(C) = \beta _F(C)$ for every circuit $C$ of
$D$.  \medskip

{\em (C)} \ There is an integer-valued potential $\pi :V\rightarrow
{\bf Z} $ for which $\underline \chi _F=\Delta _\pi $. \eL \medskip

\Proof (A)$\rightarrow $(B) \ Let $B$ be a dicut defined by a subset
$Z$ of nodes for which $\delta (Z)=0$ and $C$ a circuit.  If we go
around $C$ clockwise, and a node $v\in V-Z$ follows a node $u\in Z$,
then $vu$ is an edge of $D$, while if a node $y\in Z$ follows a node
$x\in Z$, then $xy$ is an edge of $D$.  Therefore the edges in $C\cap
B$ are alternately forward and backward edges of $C$ and hence
$\varphi _B(C) = \beta _B(C)$.  Consequently, $\varphi _F(C) = \beta
_F(C)$ holds if $F$ is the union of disjoint dicuts.

(B)$\rightarrow $(C) Let $x:=\underline \chi _F$.  Add the opposite
edge $e'$ of each edge $e$ of $D$ and define $x(e'):=-x(e)$.  Then (B)
implies that $x$ is conservative on the enlarged digraph.  By Gallai's
lemma, there is an integer-valued feasible potential $\pi $. For every
edge $e=uv\in A$ and for its opposite edge $e'=vu$, we have $\pi
(v)-\pi (u)\leq x(e)$ and $\pi (u) - \pi (v)\leq x(e')=-x(e)$ from
which $\pi (v)-\pi (u)=x(e)$, and hence $\underline \chi _F=\Delta
_\pi $.

(C)$\rightarrow $(A) Let $\pi :V\rightarrow {\bf Z} $ be a potential
for which $\underline \chi _F=\Delta _\pi $. We may assume that $D$ is
connected and also that the smallest value of $\pi $ is zero.  Let
$0=p_0<p_1<\dots <p_q$ denote the distinct values of $\pi $ and let
$Z_i:=\{v:  \pi (v)\geq p_i\}$ for $i=1,\dots ,q$.  No edge $uv$ can
leave $Z_i$, for otherwise $\pi (v)-\pi (u) \leq -1$ but $\Delta _\pi
$ is $(0,1)$-valued.  Let $B_i$ denote the set of edges entering
$Z_i$.  Since $\Delta _\pi $ is $(0,1)$-valued and $D$ is connected,
it follows that $p_i=i$.  We claim that $F=\cup B_i$.  Indeed, if
$e=uv\in F$, then $\pi (v)-\pi (u)=1$ and hence $e$ belongs to $B_i$
where $i=\pi (v)$ while if $e\in A-F$, then $\pi (v)-\pi (u)=0$ and
$e$ does not belong to any $B_i$.  $\bullet $

\medskip

Let $F\subseteq A$ be a subset of edges of $D=(V,A)$.  We say that a
cut $B$ of $D$ is {\bf $F$-clean} if every edge of $B$ in one
direction belongs to $F$ and every edge of $B$ in the other direction
belongs to $A-F$.

\Claim \label{Fclean} Let $F'$ be the symmetric difference of a subset
$F\subseteq A$ and an $F$-clean cut $B$.  Then $\varphi _F(C)=\varphi
_{F'}(C)$ for every circuit $C$ of $D$.  \eCl

\Proof If we go around $C$, then the edges of $C\cap B$ are
alternately forward and backward edges of $C$.  Since $B$ is
$F$-clean, $\varphi _F(C)=\varphi _{F'}(C)$ follows.  $\bullet $

\medskip

Two orientations $D$ and $D'$ of an undirected graph are called {\bf
dicut equivalent} if $D'$ may be obtained from $D$ by reorienting a
set of disjoint dicuts of $D'$.  Obviously, in this case $D$ can also
be obtained from $D'$ by reorienting disjoint dicuts of $D'$.  The
next lemma shows, among others, that dicut equivalence is an
equivalence relationship.

\Lemma \label{reor} Let $D=(V,A)$ and $D'=(V,A')$ be two orientations
of an undirected graph $G$.  The following are equivalent.  \medskip

{\em (A1)} \ $D$ and $D'$ are dicut equivalent.  \medskip

{\em (A2)} \ $D'$ can be obtained from $D$ by a sequence of dicut
reorientations where each time a dicut of the current member of the
sequence is reoriented.  There is a sequence where the number of
dicuts to be reoriented is at most $n-1$.  \medskip

{\em (A3)} \ $D'$ can be obtained from $D$ by a sequence of
reorienting current source-nodes.  {\em (Reorienting a source-node $v$
means that we reorient all edges leaving $v$.)} There is a sequence
where the number of reorientations is at most $(n-1)\sp 2$.  \eL

\Proof (A1)$\rightarrow $(A2) is immediate form the definition.

(A2)$\rightarrow $(A1) Suppose that $D'$ arises from $D$ as described
in (A2) and let $D_0=D,$ $\ D_1,\dots ,D_q=D'$ be a sequence of digraphs
in which each $D_{i}$ arises from $D_{i-1}$ by reorienting a dicut
$B_{i-1}$ of $D_{i-1}$.  Let $F'$ denote the subset of those edges of
$D$ which are reversed in $D'$.  We are going to show that $F'$ is the
union of disjoint dicuts of $D$.  Let $F$ denote the subset of those
edges of $D$ which are reversed in $D_{q-1}$.  By induction, $F$ is
the union of disjoint dicuts of $D$.  Let $B$ denote the cut of $D$
corresponding to the dicut $B_{q-1}$ of $D_{q-1}$.  Then $B$ is
$F$-clean and $F'$ is the symmetric difference of $B$ and $F$.  By
Claim \ref{Fclean}, $\varphi _F(C)=\varphi _{F'}(C)$ holds for every
circuit of $D$.  This implies that $\beta _F(C)=\beta _{F'}(C)$ and,
by Lemma \ref{dicuts}, $F'$ is the union of disjoint dicuts of $D$.

(A2)$\rightarrow $(A3) It suffices to show that the reorientation of a
single dicut can be obtained by a sequence of reorientations of
current source-nodes.  To this end, let $Z$ be a subset of nodes so
that no edge enters $Z$, that is, the set $B$ of edges leaving $Z$ is
a dicut.  There is a topological ordering $\{v_1,v_2,\dots ,v_n\}$ so
that the nodes of $Z$ precede the nodes outside of $Z$, that is,
$Z=\{v_1,\dots ,v_j\}$ where $j=\vert Z\vert $. Reorient first the
source-node $v_1$.  Then $v_{2}$ becomes a source-node.  Reorient now
$v_{2}$ and continue in this way until the current source-node $v_{j}$
gets reoriented.  Since each edge induced by $Z$ are reoriented
exactly twice while the edges leaving $Z$ are reoriented exactly once,
this sequence of reorientations of current source-nodes results in a
digraph that arises from $D$ by reorienting the dicut $B$.

(A3)$\rightarrow $(A2) is obvious.  $\bullet $ \medskip

Note that an acyiclic tournament $D$ on $n$ nodes shows that the bound
$(n-1)\sp 2$ in Property (C) is sharp when $D'$ is the reverse of $D$.

The property of dicut equivalence is in {\bf NP} in the sense that for
two orientations of $G$ it can be certified by exhibiting the disjoint
dicuts.  The next result provides a {\bf co-NP} characterization.

\THEOREM \label{clar1.ekvi1} Two orientations $D=(V,A)$ and
$D'=(V,A')$ of an undirected graph $G$ are dicut equivalent if and
only if

\eq \hbox{ $\varphi (C)=\varphi (C')$ for every circuit $C$ of $D$ }\
\eeq where $C'$ denotes the circuit of $D'$ corresponding to $C$.  \eT

\Proof Since the reorientation of a dicut does not change the number
of forward edges of a circuit, $\varphi (C)=\varphi (C')$ holds if $D$
and $D'$ are dicut equivalent.

Conversely, suppose that $\varphi (C)=\varphi (C')$ for every pair of
corresponding circuits.  Let $F$ denote the set of edges of $D$ that
are oppositely oriented in $D'$.  Then $\varphi _F(C) = \beta _F(C)$
for every circuit $C$ of $D$ and Lemma \ref{dicuts} implies that $F$
is the union of disjoint dicuts.  Therefore $D'$ and $D$ are dicut
equivalent.  $\bullet $ \medskip

A subset $L$ of edges of a digraph $D$ is called {\bf circuit-flat} or
just {\bf flat} if every cyclic edge of $D$ belongs to a di-circuit of
$D$ containing exactly one element of $L$.  We say that $L\subseteq A$
is a {\bf transversal} or {\bf covering} of di-circuits of $D$ if $L$
covers all di-circuits.  A subset $F\subseteq A$ of edges of a digraph
$D$ is called a {\bf flat covering} or {\bf flat transversal} of
di-circuits if $F$ covers $\cal C$ and every cyclic edge of $D$
belongs to a di-circuit covered exactly once by $F$.

\Lemma [Knuth lemma, \cite{Knuth74}] \label{Knuth} Every digraph
$D=(V,A)$ admits a flat transversal of di-circuits.  $\bullet $ \eL

Knuth formulated this result only for strongly connected digraphs but
applying his version to the strong components of $D$, one obtains
immediately the lemma.  Knuth's proof is not particularly difficult
but Iwata and Matsuda \cite{Iwata-Matsuda} found an even simpler proof
based on the ear-decomposition of strong digraphs.

\section{Characterizing the $k$-union of sink-stable and $F$-stable
sets}

By the {\bf $k$-union} of sink-stable sets, we mean a subset $U$ of
nodes that can be partitioned into $k$ sink-stable sets.  $U$ is also
called {\bf $k$-sink-stable}.  A $k$-union of $F$-stable sets is
defined analogously.  In this section, we characterize these types of
sets.  We start the investigation with $k=1$ since it behaves a bit
differently from the case $k\geq 2$.

Sink-stability was introduced as an {\bf NP}-property.  The first goal
of this section is to show that sink-stability is also in co-{\bf NP}.
That is, the next result provides an easily checkable tool to certify
that a given stable set is {\em not} sink-stable.  Recall that $\eta
(C)$ denoted the minimum of the number of forward edges and the number
of backward edges of a circuit $C$.

\THEOREM \label{sinkjel} Let $D=(V,A)$ be a digraph.  A stable set
$S\subseteq V$ is sink-stable if and only if \eq \vert S\cap V(C)\vert
\leq \eta (C) \hbox{ for every circuit $C$ of $D.$ }\
\label{(sinkjel)} \eeq \eT

\Proof If $C$ is a circuit and $v\in V(C)$ is a sink node of $D$, then
one of the two edges of $C$ entering $v$ is a forward edge and the
other one is a backward edge of $C$.  Therefore the value $\eta (C)$
is as large as the number of sink nodes in $V(C)$.  Since reorienting
a dicut does not change $\eta (C)$, we conclude that every circuit $C$
can contain at most $\eta (C)$ elements of a sink-stable set, that is,
\eref{(sinkjel)} is necessary.

To see sufficiency, assume the truth of \eref{(sinkjel)}.  Let $s\in
S$ be an element of $S$.  We can assume by induction that the elements
of $S-s$ are all sink nodes.  If no edge enters $s$, then the edges
leaving $s$ form a dicut $B$.  By reorienting $B$, the whole $S$
becomes a sink set and we are done.

Therefore we can assume that at least one edge enters $s$.  Let $T$
denote the set of nodes $u$ for which $us\in A$.  Let $D'$ denote an
auxiliary digraph arising from $D$ in such a way that we add to $D$
the opposite of all edges of $D$ entering an element of $S-s$.  Let
$Z\subseteq V$ denote the set of nodes reachable in $D'$ from $s$.
There are two cases.  \medskip

\noindent {\bf Case 1} \ $Z\cap T\not =\emptyset $, that is, $D'$
includes a directed path $P$ from $s$ to a node $t$ in $T$.  Now
$P+ts$ is a di-circuit of $D'$.  This di-circuit determines a circuit
$C$ of $D$.  If we go around $C$ in the direction of the edge $st$,
then there are exactly $\vert (S-s)\cap V(C)\vert $ oppositely
oriented edges of $C$ and there are at least $\vert (S-s)\cap
V(C)\vert $ edges in the direction of $st$.  Since $s$ belongs to $C$,
we conclude that $\vert S\cap V(C)\vert \geq \eta (C)+1$,
contradicting \eref{(sinkjel)}.

\medskip

\noindent {\bf Case 2} \ $Z\cap T=\emptyset $. Since no edge of $D'$
leaves $Z$, no edge of $D$ can leave $Z$ either.  In addition, no edge
entering $S-s$ can enter $Z$, since the opposite of such an edge
leaves $Z$ and belongs to $D'$.  Therefore the set of edges entering
$Z$ is a dicut $B$ of $D$.  By reorienting $B$, every element of $S-s$
remains a sink node.  Furthermore, $s$ becomes a source-node since the
head of each edge leaving $s$ is in $Z$ while the tail of each edge
entering $s$ is not in $Z$.  Finally, by reorienting the edges leaving
the source-node $s$, $s$ also becomes a sink node, that is, the whole
$S$ will be a sink set.  $\bullet $

\medskip Note that the proof can easily be turned to a polynomial
algorithm that either finds a circuit $C$ violating \eref{(sinkjel)}
or transforms $S$ into a sink set by reorienting a (polynomially long)
sequence of (current) dicuts, showing in this way that $S$ is a
sink-stable set.

\medskip

How can one characterize $k$-sink-stable sets when $k\geq 2$?  Before
answering this question, we recall a pretty theorem of Minty
\cite{Minty62}.  By the chromatic number $\chi(D)$ of a digraph $D$,
we simply mean the chromatic number of the underlying undirected
graph.  Minty provided an interesting upper bound for $\chi(D)$.

\THEOREM [Minty] \label{Minty} Let $D=(V,A)$ be a digraph and $k\geq
2$ an integer.  If \eq \hbox{ $ \vert C\vert \leq k \eta (C)$ for
every circuit $C$ of $D$,}\ \label{(Mintyfelt)} \eeq then $\chi(D)\leq
k$, that is, the node-set of $D$ can be partitioned into $k$ stable
sets.  $\bullet $ \eT

The theorem shows that \eref{(Mintyfelt)} is a sufficient condition
for $k$-colourability.  As a sharpening, we prove that
\eref{(Mintyfelt)} is actually a necessary and sufficient condition
for the existence of a partition of $V$ into $k$ sink-stable sets.  In
fact, we prove a bit more.

\THEOREM \label{ksink} Let $D=(V,A)$ be a digraph and $k\geq 2$ an
integer.  A subset $S\subseteq V$ is $k$-sink-stable if and only if
\eq \vert S\cap V(C)\vert \leq k\eta (C) \hbox{ for every circuit $C$
of $D$.  }\ \label{(ksink1)} \eeq \eT

\Proof We have already observed in Theorem \ref{sinkjel} that a
circuit $C$ can contain at most $\eta (K)$ elements of a sink-stable
set from which the necessity of \eref{(ksink1)} follows.

To see sufficiency, consider the digraph $D\sp *=(V,A \cup A')$
arising from $D$ by adding the reverse of every edge of $D$.  Define a
cost function $c$ on $A\cup A'$ as follows.  For an edge $a$ of $D$,
let $c(a)=k$ and for the reverse $a'$ of $a$ let $c(a')=0$.  For a
two-element di-circuit $K$ consisting of edges $a$ and $a'$ we have
$\vert S\cap V(K)\vert \leq 2\widetilde c(K)$ holds since $k\geq 2$.
Hence \eref{(ksink1)} is equivalent to the following condition.

\eq \vert S\cap V(K)\vert \leq \widetilde c(K) \hbox{ for every
di-circuit $K$ of $D\sp *$.}\ \label{(ksink2)} \eeq

Revise now $c$ in such a way that $c(e)$ is reduced by 1 for every
edge of $D\sp *$ for which the head is in $S$.  Let $c\sp *$ denote
the resulting cost function.  Observe that the $\widetilde c\sp
*$-cost of a di-circuit $K$ of $D\sp *$ is equal to $\widetilde c(K)$
minus the number of edges of $K$ having their head in $S$, that is,
$\widetilde c\sp *(K)=\widetilde c(K) - \vert S\cap V(K)\vert $.
Therefore \eref{(ksink2)} is equivalent to requiring that $c\sp *$ is
conservative.

By Lemma \ref{konzerv} there is an integer-valued $c\sp *$-feasible
potential $\pi $. Since $\pi $ can be translated by a constant, we can
assume that the smallest value of $\pi $ is 0. Let $M$ denote the
maximum value of $\pi $ and consider the following sets for $i=0,\dots
,M$.

$$P_i:=\{v:  \pi (v)=i\} \ \ \hbox{ and}\ \ \ U_i:=P_0\cup P_1\cup
\cdots \cup P_i.$$

Moreover, define for $j=0,\dots ,k-1$ the following sets.

$$V_j:=\{v:  \pi (v)\equiv j \hbox{ mod}\ k\}\ \ \hbox{ and}\ \ \
S_j:=V_j\cap S.$$

For each $uv\in A$, we have $\pi (v)\geq \pi (u)$ since $c\sp
*(vu)\leq 0$ from which $\pi (u)-\pi (v)\leq c\sp *(vu)\leq 0$.
Therefore no edge of $D$ enters any $U_i$, that is, the set $B_i$ of
edges of $D$ leaving $U_i$ is a dicut of $D$.  Obviously, the sets
$S_j$ partition $S$.  We are going to prove that each $S_j$ is a
sink-stable set from which the theorem will follow.  To this end,
consider the dicuts $B_j, B_{j+k}, B_{j+2k},\dots $. These are
disjoint since $\pi (v)-\pi (u) \leq c\sp *(uv)\leq k$ holds for each
edge $uv\in A$.

Let $z\in S_j$.  For any edge $uz\in A$ entering $z$, we have $\pi
(z)-\pi (u)\leq c\sp *(uz) = k-1$ and hence $uz$ is not in any of the
dicuts $B_j, B_{j+k}, B_{j+2k},\dots $. For any edge $zv\in A$ leaving
$z$, we have $c\sp *(vz)=-1$ from which $\pi (z) - \pi (v) \leq c\sp
*(vz)=-1$ and hence $\pi (v)- \pi (z) \geq 1$.  Therefore $zv$ belongs
to one of the dicuts $B_j, B_{j+k}, B_{j+2k},\dots $. Consequently,
each node $z\in S_j$ is a sink node in $D_B$ where $B$ is the union of
the dicuts $B_j, B_{j+k}, B_{j+2k},\dots $ and $D_B$ denotes the
digraph arising from $D$ by reversing $B$.  $\bullet $

\medskip

It is known that there is a polynomial time algorithm for an arbitrary
cost function $c$ that either finds a $c$-feasible potential or finds
a negative di-circuit.  Therefore the proof of Theorem \ref{ksink}
above gives rise to an algorithm that either finds a partition of $S$
into $k$ sink-stable sets or finds a circuit $C$ of $D$ violating
\eref{(ksink1)}.  \medskip

{\bf Remark} \ It is useful to observe that for $k=1$ the statement in
Theorem \ref{ksink} fails to hold:  in a digraph consisting of two
nodes and a single edge, \eref{(ksink1)} holds automatically since
there is no circuit at all but $V$ is not a stable set.  This is why
we assumed a priori in Theorem \ref{sinkjel} that $S$ is a stable set.
We also remark that the proof technique of Theorem \ref{ksink} can be
used for $k=1$, as well, to obtain an alternative proof for the
non-trivial direction of Theorem \ref{sinkjel} since in the latter
case the stability of $S$ is part of the assumption.  \medskip

Next, we describe a characterization for the union of $k$ $F$-stable
sets.  Since $F$-stability was defined through sink-stability, it is a
straightforward task to translate the theorems above on characterizing
$k$-unions of sink-stable sets to those on characterizing $k$-unions
of $F$-stable sets.  The resulting necessary and sufficient condition
is a certain inequality required to hold for {\em all} circuits of the
digraph.  In the applications, however, the digraph in question is
strongly connected, and in this case it turns out that it suffices to
require an inequality to hold only for every di-circuit.

For this simplification, we shall need a lemma.  Let $F\subseteq A\sp
*$ be a subset of edges of a digraph $D\sp *=(V,A\sp *)$.  For a
closed walk $W$, let $$\sigma _F(W):= \varphi _F(W) + \beta _{A\sp
*-F}(W)$$ where $\varphi _F(W)$ and $\beta _{A\sp *-F}(W)$ denote the
number of forward $F$-edges and the backward $(A\sp *-F)$-edges,
respectively.  Therefore $\sigma _F(W)= \vert F\cap W\vert $ for a
one-way walk $W$ and, in particular, for a di-circuit.  Recall that
$F$ was called flat if every cyclic edge of $D\sp *$ belongs to a
di-circuit covered by exactly one element of $F$.

\Lemma \label{seta} If $F\subseteq A\sp *$ is a flat subset of a
digraph $D\sp *=(V,A\sp *)$, then the node-set $V(C)$ of every circuit
$C$ of $D\sp *$ consisting of cyclic edges can be covered by
di-circuits $K_1,\dots ,K_q$ for which $\sum _i \vert F\cap K_i\vert =
\sum _i \sigma _F(K_i) \leq \sigma _F(C)$.  \eL

\Proof It is more convenient to prove the more general statement
asserting that the node-set of a closed walk $W$ consisting of cyclic
edges of $D\sp *$ can be covered by di-circuits $K_1,\dots ,K_q$ for
which $\sum _i \vert F\cap K_i\vert = \sum _i \sigma _F(K_i) \leq
\sigma _F(W)$.  We use induction on the number of backward edges of
$W$.  If this number is zero, that is if $W$ is a one-way walk, then
by traversing all the edges of $W$, we obtain di-circuits $K_1,\dots
,K_q$ for which $\sum _i\vert F\cap K_i\vert = \sum _i \sigma _F(K_i)
= \sigma _F(C)$.

Suppose now that $e=uv\in A\sp *$ is a backward edge of $W$.  By the
hypothesis, $e$ belongs to a di-circuit $K$ containing one $F$-edge,
that is, there is a directed path $P$ from $v$ to $u$.  Replace $e$ in
the walk $W$ by $P$ and let $W'$ denote the closed walk obtained in
this way.  Obviously $V(W)\subseteq V(W')$ and $W'$ has one less
backward edge than $W$.  We claim that $\sigma _F(W') \leq \sigma
_F(W)$.  Indeed, if $e\in F$, then $e$ contributes to $\sigma _F(W)$
by zero and $P$ contains no $F$-edge from which $\sigma _F(W') =
\sigma _F(W)$.  If $e\in A\sp *-F$, then $e$ contributes to $\sigma
_F(W)$ by 1. Furthermore, since $P$ contains one $F$-edge, we confer
$\sigma _F(W') = \sigma _F(W)$ from which the lemma follows by
induction.  $\bullet $

\THEOREM \label{Fsinkjel} Let $F\subseteq A\sp *$ be a flat subset of
edges of a strongly connected digraph $D\sp *=(V,A\sp *)$ and let
$k\geq 1$ be an integer.  A subset $S\subseteq V$ is the union of $k$
$F$-stable sets if and only if \eq \vert S\cap V(K)\vert \leq k\vert
F\cap K\vert \hbox{ for every di-circuit $K$ of $D\sp *$.  }\
\label{(Fsinkjel)} \eeq In particular, if $F$ is a flat transversal of
di-circuits, then the minimum number of $F$-stable sets covering $S$
is equal to \eq \max \left \{ \left \lceil { \vert S\cap V(K)\vert
\over \vert F\cap K\vert } \right \rceil :  K \hbox{ a di-circuit of}\
D\sp * \right \}.  \label{(Fsinkjel2)} \eeq \eT

\Proof Let $D'=D\sp *_F$ denote the digraph arising from $D\sp *$ by
reversing $F$.  Suppose that $S$ is the $k$-union of $F$-stable sets,
that is, $S$ is $k$-sink-stable in $D'$.  Let $K$ be a di-circuit of
$D\sp *$ and let $K'$ denote the corresponding circuit of $D'$.  Then
$\vert S\cap V(K)\vert \leq k \eta (K')\leq k\vert F\cap K\vert $,
from which the necessity of \eref{(Fsinkjel)} follows.

Suppose now that \eref{(Fsinkjel)} holds.  For every set $X$ of edges
of $D\sp *$, the corresponding set in $D'$ will be denoted by $X'$.

\Claim \label{Fsinkjelclaim} $\vert S\cap V(C')\vert \leq k\eta (C')$
for every circuit $C'$ of $D'$.  \eCl

\Proof We may assume that $\eta (C') = \beta (C') \leq \varphi (C')$.
Let $C$ be the circuit of $D\sp *$ corresponding to $C'$.  By applying
Lemma \ref{seta} to $D\sp *$ and to $C$, we obtain that $V(C)=V(C')$
can be covered by di-circuits $K_1,\dots ,K_q$ of $D\sp *$ for which
$\sum _i\vert F\cap K_i\vert = \sum _i \sigma _{F}(K_i) \leq \sigma
_F(C)$.  By applying \eref{(Fsinkjel)} to di-circuits $K_i$, we see
that $\vert S\cap V(K_i)\vert \leq k\vert F\cap K_i\vert $. Hence
$\vert S\cap V(C')\vert = \vert S\cap V(C)\vert \leq \sum _i\vert
S\cap V(K_i)\vert \leq \sum _i k\vert F\cap K_i\vert \leq k \sigma
_F(C) = k [\varphi _F(C) + \beta _{A\sp *-F}(C)] = k\beta (C') = k\eta
(C')$, as required.  $\bullet $

\medskip If $k=1$, then every edge of $D\sp *$ belongs to a di-circuit
$K$ for which $\vert F\cap K\vert =1$ and hence $S$ is stable by
\eref{(Fsinkjel)}.  Theorem \ref{sinkjel}, when applied to $D'$,
implies that $S$ is sink-stable in $D'$, that is, $S$ is $F$-stable in
$D\sp *$.  If $k\geq 2$, then we can apply Theorem \ref{ksink} to
$D'$.  By Claim \ref{Fsinkjelclaim} above, $S$ is $k$-sink-stable in
$D'=D\sp *_F$, showing that $S$ is the $k$-union of $F$-stable sets in
$D\sp *$.  $\bullet $ $\bullet $

\medskip

Note that, unlike the corresponding situation with sink-stable sets
where we formulated the characterization of $k$-unions of sink-stable
sets separately for $k=1$ and for $k\geq 2$, in the formulation of
Theorem \ref{Fsinkjel} these cases are not separated.  It is the proof
of the theorem where the two cases were handled separately.

\Corollary [Bondy, \cite{Bondy76}] The chromatic number $\chi (D)$ of
a strongly connected digraph $D\sp *=(V,A\sp *)$ is at most the length
of the longest di-circuit of $D\sp *$.  \eCo

\Proof Let $F$ be a flat transversal of di-circuits ensured by Lemma
\ref{Knuth}.  By applying Theorem \ref{Fsinkjel} to $S:=V$, we confer
that $\chi(D) \leq $ the minimum number of $F$-stable sets covering $V
= \max \{ \lceil { \vert K\vert \over \vert F\cap K\vert } \rceil :
K$ a di-circuit of $D\sp * \} \leq \max \{ \vert K\vert :  K$ a
di-circuit of $D\sp *\}$, as required.  $\bullet $

\medskip In the last section, we will point out a link between
$F$-stable sets and so-called cyclic stable sets introduced by Bessy
and Thomass\'e \cite{Bessy-Thomasse}.  It was their paper that first
showed how Bondy's theorem follows from results on cyclic stable sets.

\medskip

{\bf Remark} \ We showed above how Theorems \ref{sinkjel} and
\ref{ksink} gave rise to Theorem \ref{Fsinkjel}.  But the reverse
derivation is also possible.  To obtain, for example, the non-trivial
sufficiency part of Theorem \ref{ksink}, suppose that \eref{(ksink1)}
holds.  Then every di-circuit of $D$ is disjoint from $S$.  Let $F$
denote the set of reverse edges of $D$ and consider the strong digraph
$D\sp +=(V, A+F)$.  Then $F$ is clearly flat since each original edge
of $D$ and its reverse edge form a di-circuit covered once by $F$.

Let $K$ be a di-circuit of $D\sp +$.  We are going to show that
\eref{(Fsinkjel)} holds.  If $K$ consists of edges of $D$, that is, if
$F\cap K=\emptyset $, then $S$ is disjoint from $V(K)$ and hence $
\vert S\cap V(K)\vert = 0 \leq k \vert F\cap K\vert $. Suppose now
that $F\cap K\not =\emptyset $. If $\vert K\vert =2$, then $\vert
K\vert =2 \leq k \vert F\cap K\vert $. Finally, if $\vert K\vert \geq
3$, then $K$ determines a circuit $C$ of $D$ so that the backward
edges of $C$ correspond to the elements of $F\cap K$.  Since $\vert
S\cap V(C)\vert \leq k \eta (C)$ by \eref{(ksink1)}, we obtain that
$\vert S\cap V(K)\vert \leq k \vert F\cap K\vert $.

By applying Theorem \ref{Fsinkjel} to $D\sp *$ in place of $D$, we
obtain that $S$ is the union of $k$ $F$-stable sets.  By definition a
set $Z$ is $F$-stable if $S$ is sink-stable in $D\sp +_F$.  Since
$D\sp +_F$ arises from $D$ by duplicating each edge in parallel, $Z$
is sink-stable in $D$, as well, and hence $S$ is indeed the union of
$k$ sink-stable sets.

Theorem \ref{sinkjel} follows from the special case $k=1$ of Theorem
\ref{Fsinkjel} in an analogous (and even simpler) way.  \medskip

\section{Optimal sink-stable and $F$-stable sets}

\medskip

Our next goal is to investigate sink-stable and $F$-stable sets of
largest cardinality and, more generally, of maximum weight.  The main
device to obtain min-max theorems for these parameters is a result of
Gallai.

\medskip

\noindent {\bf A theorem of Gallai} \medskip

Let $D=(V,A)$ be a digraph and $c:A\rightarrow {\bf Z}_+$ a
non-negative integer-valued function.  The {\bf $c$-value} of a
circuit $C$ is the sum of the $c$-values of the edges of $C$, that is,
$\widetilde c(C)$.  We say that a multiset of nodes is {\bf
$c$-independent} if it contains at most $\widetilde c(K)$ nodes of
every di-circuit $K$ of $D$.  A multiset can be identified with a
non-negative integral vector $x:V\rightarrow {\bf Z}_+$ and then the
$c$-independence of $x$ means that $\widetilde x(V(K)) \leq \widetilde
c(K)$ for every di-circuit $K$.

Let $w:V\rightarrow {\bf Z}_+$ be a weight function.  For a function
$y\geq 0$ defined on the set of di-circuits of $D$, we say that $y$
{\bf covers} $w$ if

\eq \hbox{ $\sum [y(K):  v\in V(K), K$ a di-circuit$]\geq w(v)$ for
every $v\in V$.}\ \eeq

A circulation $z\geq 0$ is said to cover $w$ if $\varrho _z(v)\geq
w(v)$ holds for every node $v\in V$.  The following lemma describes a
simple and well-known relationship between circulations and families
of di-circuits covering $w$.

\Lemma \label{gallema} If $y\geq 0$ is a function on the set of
di-circuits covering $w$, then the function $z:A\rightarrow {\bf Z}_+$
defined by $z(e):  = \sum [y(K):  K$ a di-circuit and $e\in K]$ is a
non-negative circulation covering $w$ for which $cz= \sum
[y(K)\widetilde c(K):  K$ a di-circuit$]$.  Furthermore, if $y$ is
integer-valued, then so is $z$.  Conversely, a circulation $z\geq 0$
covering $w$ can be expressed as a non-negative linear combination of
di-circuits, and if $z= \sum y(K)\underline \chi (K)$ is such an
expression, then $y$ covers $w$ and $cz= \sum [y(K)\widetilde c(K):
K$ a di-circuit$]$.  Furthermore, if $z$ is integer-valued, then $y$
can also be chosen integer-valued.  $\bullet $ \eL

The following result of Gallai \cite{Gallai58} appeared in 1958.  We
cite it in its original form because the literature does not seem to
know about it.  \medskip

{\bf (3.2.7) SATZ.} \ \ {\em Ist $\Gamma $ endlich und gilt
$\psi[k]\geq 0$ f\"ur jeden positiven Kreis $k$ gibt ist ferner zu
jedem Punkt $X$ mit $\varphi (X)>0$ einen positiven Kreis, der $X$
enth\"alt, so ist das Minimum der $\psi $-Werte der punktf\"ullende
positiven Kreissysteme gleich dem Maximum der $\varphi $-Werte der
kreisaufnehmbaren Punktsysteme.}

\medskip

That is:  \ {\em If $\Gamma $ is finite and $\psi [k]\geq 0$ holds for
every positive circuit, and if, furthermore, every point $X$ with
$\varphi (X)>0$ is included in a positive circuit, then the minimum
$\psi $-value of point-covering positive circuit-systems is equal to
the maximum $\varphi $-value of circuit-independent point-systems.}
\medskip

Here $\psi$ and $\varphi $ are integer-valued functions on the
edge-set and on the node-set, respectively, of the digraph $\Gamma $,
and a positive circuit means a di-circuit.  In the present context, we
use functions $c$ and $w$ in place of $\psi$ and $\varphi $,
respectively, and the theorem can be formulated as follows.

\THEOREM [Gallai, Theorem (3.2.7) in \cite{Gallai58}] \label{Gallaiw}
Let $c:A\rightarrow {\bf Z}_+$ and $w:V\rightarrow {\bf Z}_+$ be
non-negative functions on the edge-set and on the node-set,
respectively, of a digraph $D=(V,A)$, and assume that for each $v\in
V$ with $w(v)>0$ that $v$ belongs to a di-circuit.  Then the minimum
total sum of $c$-values of a system of di-circuits covering $w$ is
equal to the maximum $w$-weight of a $c$-independent multiset of nodes
of $D$, or more formally, to \eq \hbox{ $\max \{wz:  z\in {\bf Z}_+\sp
V, \ z \ c$-independent$\}.$ }\ \label{(clar3.max)} \eeq \eT

Note that in the original version cited above only the
conservativeness of $c \ (=\psi)$ was assumed and not its
non-negativity.  But for a conservative $c$ there is a feasible
potential $\pi $ and then the cost function $c_\pi $ defined by $c_\pi
(uv)= c(uv) - \pi (v)+\pi (u)$ is non-negative for which $\widetilde
c(K)=\widetilde c_\pi (K)$ holds for every di-circuit $K$.  That is,
the theorem for conservative $c$ follows from its special case for
non-negative $c$.

In the special case of Theorem \ref{Gallaiw}, when $w$ is
$(0,1)$-valued, that is, when $w:=\underline \chi _U$ for a subset
$U\subseteq V$ of nodes, we have the following.

\Corollary \label{GallaiU} Let $U\subseteq V$ be a specified subset of
nodes belonging to some di-circuits.  The minimum total $c$-value of
di-circuits covering $U$ is equal to the maximum number of {\em
(not-necessarily distinct)} $c$-independent elements of $U$.  $\bullet
$ \eCo

For completeness, we outline a proof of Theorem \ref{Gallaiw}.  Let
$Q$ denote the di-circuits versus nodes incidence matrix of a digraph
$D$.  That is, $Q$ is a $(0,1)$-matrix with rows corresponding to the
di-circuits and columns corresponding to nodes.  An entry
corresponding to a di-circuit $C$ and a node $v$ is 1 or zero
according to whether $v$ is in $V(C)$ or not.  A fundamental theorem
of Edmonds and Giles \cite{Edmonds-Giles} states that a polyhedron $R$
is integral provided that $R$ is described by a totally dually
integral system and both the constraint matrix and the bounding vector
is integral.  Combined with the l.p. duality theorem, the following
result is equivalent to Theorem \ref{Gallaiw}.

\THEOREM \label{Gallaitdi} Let $Q$ be the di-circuits versus node
incidence matrix of a digraph $D$.  Let $\widetilde c$ denote a vector
the components of which correspond to the rows of $Q$ {\em (that is,
to the di-circuits of $D$)} and the value of a component corresponding
to a di-circuit $K$ is $\widetilde c(K)$.  Then the linear system \eq
\{Qx \leq \widetilde c, \ x\geq 0\} \label{(TDI)} \eeq is totally dual
integral.  \eT

\Proof Let $w$ be an integer-valued function on $V$.  Consider the
following linear program:

\eq \hbox{ $ \min \{\sum [y(K)\widetilde c(K):  K$ a di-circuit$] :
yQ \geq w, \ y\geq 0\}.$ }\ \label{(dual)} \eeq

What we have to show is that this program has an integral-valued
optimum if it has an optimum at all.  We may assume that $w$ is
non-negative.  We may also assume that each $v\in V$ with $w(v)>0$
belongs to a di-circuit for otherwise \eref{(dual)} has no feasible
solution at all.

By Lemma \ref{gallema}, it suffices to show that the linear system
$\min \{cz:  z\geq 0$ a circulation covering $w\}$ has an
integer-valued optimum.  But this follows from the integrality of the
circulation polyhedron by appplying the standard node-duplicating
technique.  Indeed, replace each node $v$ of $D$ by nodes $v'$ and
$v''$, replace each edge $uv\in A$ by a new edge $u'v''$ (with lower
capacity $0$ and cost $c(uv)$), and finally add a new edge $v''v'$
(with lower capacity $w(v)$ and cost 0) for every original node $v\in
V$.  In the resulting digraph $D'$, a feasible circulation $z'$
defines a non-negative circulation $z$ of $D$ which covers $w$ (that
is, $\varrho _z(v)\geq w(v)$ for $v\in V$) and $c'z'=cz$.  $\bullet $
\medskip

With some work, this proof can be used to turn a min-cost circulation
algorithm to one that computes the optima in Theorem \ref{Gallaiw} in
polynomial time.  \medskip

\noindent {\bf On a solved conjecture of Gallai}

\medskip

Before turning to $F$-stable sets, we make a detour and show how
Gallai's theorem from 1958 and and Knuth's lemma from 1974 imply
immediately the following conjecture of Gallai \cite{Gallai63} that
was first proved by Bessy and Thomass\'e \cite{Bessy-Thomasse}.
Recall that $\alpha (D)$ denotes the stability number of $D$ while
$\gamma (D)$ denotes the minimum number of di-circuits fo $D$ covering
$V$.

\THEOREM [Bessy and Thomass\'e] \label{Gallaisejt} Let $D=(V,A)$ be a
strongly connected digraph with at least two nodes.  Then $\gamma (F)
\leq \alpha (D)$, that is, $V$ can be covered by $\alpha (D)$
di-circuits.\eT

\Proof By Lemma \ref{Knuth} $D$ has a flat covering $F$ of
di-circuits.  By applying Corollary \ref{GallaiU} to $U:=V$ and to
$c:=\underline \chi _F$, we obtain that the minimum total $c$-weight
$\gamma _F$ of di-circuits covering $V$ is equal to the maximum number
$\alpha _F$ of $c$-independent elements of $V$.  Since $F$ covers
every di-circuit $K$, the $c$-weight of $K$ is at least 1 and hence
$\gamma \leq \gamma _F$.  Since $F$ is a flat covering of di-circuits,
every node belongs to a di-circuit of $c$-cost 1. Therefore a
$c$-independent multiset $S$ is actually a set.  The $c$-independence
also implies that $S$ is a stable set.  Hence $\gamma (D)\leq \gamma
_F = \vert S\vert \leq \alpha (D)$.  $\bullet $ \medskip

Note that the same argument shows that the following extension also
holds.

\THEOREM Let $D=(V,A)$ be a strongly connected digraph and $U$ a
subset of nodes of $D$.  Then $U$ can be covered by $\alpha _U$
di-circuits where $\alpha _U$ denotes the maximum cardinality of a
stable subset of $U$.  $\bullet $ \eT

We note that Cameron and Edmonds \cite{Cameron-Edmonds92} (not knowing
of the paper of Gallai \cite{Gallai58}) proved an extension of Theorem
\ref{Gallaitdi} asserting that the linear system $\{Qx\leq \widetilde
c\}$ is actually box-TDI (see, Theorem \ref{CE} below).  Based on
this, they derived Theorem \ref{Gallaisejt} in
\cite{Cameron-Edmonds08}.

\medskip

\noindent {\bf Optimal $F$-stable and sink-stable sets} \medskip

Let $F$ be a flat subset of edges of a strongly connected digraph
$D\sp *=(V,A\sp *)$.  Here $F$ is not necessarily a transversal of
di-circuits.  Let $Q$ denote the di-circuit versus node incidence
matrix of $D\sp *$ and let $c_F$ be a vector whose components
correspond to the di-circuits of $D\sp *$ and $c_F(K)$ is the
$F$-value $\vert F\cap K\vert $ of $K$ for a di-circuit $K$.

\THEOREM \label{Fmax.clar3} Let $F$ be a flat subset of edges of a
strongly connected digraph $D\sp *=(V,A\sp *)$ and let $w:V\rightarrow
{\bf Z}_+$ be an integer-valued weight-function on the node-set of
$D\sp *$.  \ The maximum $w$-weight of an $F$-stable set of $D\sp *$
is equal to the minimum total $F$-value of di-circuits of $D\sp *$
covering $w$.  For a given subset $U\subseteq V$, the maximum
cardinality of an $F$-stable subset of $U$ is equal to the the minimum
total $F$-value of di-circuits of $D\sp *$ covering $U$.  \eT

\Proof Apply Gallai's theorem to $c_F$.  Since every edge of $D\sp *$
belongs to a di-circuit $K$ for which $c_F(K)=1$, a $c_F$-independent
multiset of nodes of $D\sp *$ is actually a set $S\subseteq V$, and by
Theorem \ref{Fsinkjel} $S$ is $F$-stable.  Therefore the result is a
direct consequence of Theorem \ref{Gallaiw}.  The second half of the
theorem follows by applying the the first one to $w:=\underline \chi
_U$.  $\bullet $

\THEOREM \label{sinkmax} Let $D=(V,A)$ be a digraph with no isolated
nodes and let $w:V\rightarrow {\bf Z}_+$ be an integer-valued
weight-function on the node-set of $D$.  The maximum $w$-weight of a
sink-stable set of $D$ is equal to the minimum total value of circuits
and edges of $D$ covering $w$ where the value of a circuit $C$ is
$\eta (C)$ while the value of an edge is 1. For a given subset
$U\subseteq V$, the maximum cardinality of a sink-stable subset of $U$
is equal to the minimum total value of circuits and edges of $D$
covering $U$.  \eT

\Proof For a sink-stable subset $S$, an edge can cover at most one
element of $S$.  In Theorem \ref{sinkjel} we already observed that a
circuit $C$ can cover at most $\eta (C)$ elements of a sink-stable set
from which $\max \leq \min$ follows.

The proof of the reverse direction $\max \geq \min$ can be made
separately for the components of $D$, and hence we can assume that $D$
is weakly connected.  Let $D\sp *=(V,A\cup A')$ be the digraph arising
from $D$ by adding the reverse of each edge of $D$.  Here $A'$ denotes
the set of reverse edges of $D$.  $D\sp *$ is clearly strongly
connected and $F:=A'$ is flat since each edge $e\in A$ and its reverse
$e'\in A'$ form a 2-element di-circuit covered once by $F$.

There are two types of di-circuits of $D\sp *$.  Type I is of form
$K=\{e,e'\}$ where $e=uv\in A$ and $e'=vu\in A'$, and in this case
$\widetilde c_F(K)=1$.  A Type II di-circuit $K$ arises from a circuit
$C$ of $D$ by reversing its forward edges or by reversing its backward
edges.  Therefore if $K$ is such a di-circuit of $D\sp *$, then the
reverse $\vec K$ of $K$ is also a di-circuit of $D\sp *$, and $\eta
(C) = \min \{ \widetilde c_F(K), \ \widetilde c_F(\vec K)\}$.  For
notational convenience, we will assume that $\eta (C) = \widetilde
c_F(K)$.

Let $S$ be a sink-stable set of $D$ with maximum $w$-weight.  Since
$D\sp *_F$ is a digraph that can be obtained from $D$ by doubling each
edge of $D$ in parallel, a subset $Z$ of nodes is sink-stable in $D$
if and only of $Z$ is $F$-stable in $D\sp *$.  Therefore $S$ is an
$F$-stable set in $D\sp *$ of maximum $w$-weight.

In order to prove $\max \geq \min$, we are going to show that there is
a family of circuits and edges of $D$ covering $w$ for which the total
value is $\widetilde w(S)$.  Since $S$ is a maximum $w$-weight
$F$-stable set, Theorem \ref{Fmax.clar3} implies the existence of a
family ${\cal C}\sp *$ of di-circuits of $D\sp *$ covering $w$ for
which the total $F$-value is $w(S)$.  As mentioned above, a di-circuit
in $D\sp *$ of Type I determines an edge of $D$ while a di-circuit $K$
of Type II determines a circuit $C$ of $D$ for which $\eta (C)=
c_F(K)$.  Therefore ${\cal C}\sp *$ defines a family of edges and
circuits of $D$ covering $w$ for which the total value is the total
$c_F$-value of ${\cal C}\sp *$, that is $\widetilde w(S)$, and hence
the requested direction $\max \geq \min$ follows.

The second half of the theorem follows by applying the the first one
to $w:=\underline \chi _U$.  $\bullet $

\medskip

\noindent {\bf Remark} \ One may be wondering whether the minimal
covering of $U$ in Theorem \ref{sinkmax} can perhaps be realized only
by circuits, without using edges.  The following example shows,
however, that the use of edges is anavoidable.  Let
$U:=V:=\{a,b,c,d,e\}$ and let the edges of $D$ be
$\{ab,ac,ad,eb,ec,ed\}$.  In this digraph $S=\{b,c,d\}$ is a largest
sink-stable set.  On the other hand, each circuit $C$ of $D$ has 4
edges and $\eta (C)=2$.  Therefore the total value of the best
covering of $V$ by only circuit is 4. An optimal covering consists of
a circuit with edge-set $\{ab,be,ec,ca\}$ and of an edge $ad$ with
total value 3.

\begin{figure}[h]
\begin{center}
{\includegraphics[scale=0.6,angle=270]{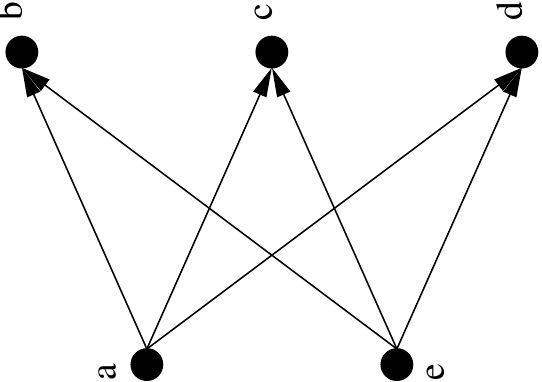}}%
\caption{\label{fig:ex_cover} The largest sink-stable set in this graph is $S =
  \{b,c,d\}$. An optimal covering consists of a circuit $C = \{ab, be, ec,
  ca\}$ and an edge $ad$, the total value of this covering is $3$.}
\end{center}
\vspace{-3ex}
\end{figure}

\section{Clar number of plane bipartite graphs}

As an application of Theorem \ref{Fmax.clar3}, we derive a recent
min-max theorem of Abeledo and Atkinson \cite{Abeledo-Atkinson1} on
the Clar number of bipartite plane graphs.  Let $G=(S,T;E)$ be a
perfectly matchable 2-connected bipartite plane graph.  The expression
plane graph means that $G$ is planar and we consider a fixed embedding
in the plane.  Note that the embedding subdivides the plane into
regions, among them exactly one is unbounded.  The bounded regions
will be referred to as faces of $G$.  Since $G$ is 2-connected, each
region is bounded by a circuit of $G$.

We call a set of faces {\bf resonant} if the faces are disjoint and
their deletion leaves a perfectly matchable graph.  Here the deletion
of a face means that we delete all the nodes of the circuit bounding
the face.  The {\bf Clar number} of $G$ is defined to be the maximum
cardinality of a resonant set of faces.  For example, if $G$ is the
graph of a cube, then the Clar number is 2, independently of the
embedding.  It is not difficult to find an example where the Clar
number does depend on the embedding.  This notion was originally
introduced in chemistry by E. Clar \cite{Clar-book} for hexagonal
graphs (where each face is a 6-circuit) to capture the behaviour of
characteristic chemical and physical properties of aromatic
benzenoids.

Before stating the result of Abeledo and Atkinson on the Clar number,
we introduce some notation that will be used in the theorem and also
in its proof.  With the bipartite graph $G=(S,T;E)$, we associate a
digraph $D=(V,A)$, where $V=S\cup T$, arising from $G$ by orienting
all edges from $S$ to $T$.  Clearly, $D$ is acyclic.  A subset of
edges of $D$ corresponding to a subset $X$ of edges of $G$ will be
denoted by $\vec X$.  For a digraph $H$ and subset $F$ of its edges,
$H_F$ denotes the digraph arising from $H$ by reversing $F$.  For a
subset $Z\subseteq V$, the set $B=\Delta (Z)$ of edges connecting $Z$
and $V-Z$ will be called a cut of $G$ determined by $Z$.  We call such
a cut of $G$ {\bf feasible} if it determines a dicut in the associated
digraph $D$.  The {\bf value} val$(B)$ of $B$ is defined to be the
absolute value of $\vert S\cap Z\vert -\vert S\cap T\vert $. It is an
easy exercise to see that val$(B)= d_M(Z)=\vert M\cap B\vert $ for an
arbitrary perfect matching $M$ of $G$.  In particular, this means for
a feasible cut $\Delta (Z)$ that $d_M(Z)$ is independent of the choice
of perfect matching $M$.

Let $G\sp *=(V\sp *,E\sp *)$ denote the planar dual of $G$.  There is
a one-to-one correspondence between the regions of $G$ and the nodes
of $G\sp *$ and also there is a one-to-one correspondence between the
edges of $G$ and the edges of $G\sp *$.  We use the convention that
for a subset $X\subseteq E$ of edges, the corresponding subset of
edges of $G\sp *$ is denoted by $X\sp *$.  It is well-known that $X$
is a circuit of $G$ if and only if $X\sp *$ is a bond of $G\sp *$.  (A
bond of a graph is a minimal cut.  A useful property is that a cut
$B=\Delta (Z)$ of a connected graph is a bond if and only if both $Z$
and $V-Z$ induce a connected subgraph, and another one is that every
cut can be partitioned into bonds.)

We also need the planar dual digraph $D\sp *=(V\sp *,A\sp *)$ of $D$.
This is an orientation of $G\sp *$ in such a way that, for a pair
$e\in E$ and $e\sp *\in E\sp *$ of corresponding edges, if the
directed edge $\vec e$ of $D$ is represented in the plane by a
vertical line segment oriented downward, then the corresponding
horizontal dual edge $\vec e\sp *$ of $D\sp *$ is oriented from right
to left.  It follows that a subset $X$ of edges of $D$ is a minimal
dicut if and only if the corresponding subset $X\sp *$ of edges of
$D\sp *$ is a di-circuit.  Similarly, a subset $X$ of edges of
$D_{\vec M}$ is a di-circuit if and only if the corresponding set
$X\sp *$ of edges of $D\sp *_{\vec M\sp *}$ is a minimal dicut of
$D\sp *_{\vec M\sp *}$.  In particular, the circuit bounding a region
of $D_{\vec M}$ is a di-circuit oriented clockwise if and only if the
corresponding node of $D\sp *_{\vec M\sp *}$ is a sink node.  Since
$D$ is acyclic, $D\sp *$ is strongly connected.  Since each node in
$T$ determines a dicut of $D$, each edge of $D$ belongs to a dicut
covered exactly once by $\vec M$. Therefore $\vec M\sp *$ is flat in
$D\sp *$.  Note, however, that $\vec M$ need not cover all dicuts of
$D$ and hence $\vec M\sp *$ is not necessarily a transversal of
di-circuits of $D\sp *$.

\THEOREM [Abeledo and Atkinson] \label{Abeledo-Atkinson} Let
$G=(S,T;E)$ be a 2-connected perfectly matchable plane bipartite
graph.  The Clar number of $G$ is equal to the minimum total value
of feasible cuts intersecting all faces of $G$.\eT

\Proof Let $M$ be an arbitrary perfect matching of $G$.  Consider the
digraph $D\sp *_{F}$ which is the planar dual digraph of $D_{\vec M}$,
where $F:=\vec M\sp *$.

\Lemma \label{equi} A set $\cal S$ of disjoint faces of $G$ is
resonant if and only if the corresponding set $S$ of nodes of $D\sp *$
is $F$-stable.  \eL

\Proof By definition, $\cal S$ is resonant if there is a perfect
matching $M'$ of $G$ so that the bounding circuit of each member of
$\cal S$ is $M'$-alternating.  This is equivalent to requiring that
these bounding circuits are directed circuits in $D_{\vec M'}$.  By
reorienting such a di-circuit if necessary, we can assume that the
bounding circuits of the members of $\cal S$ are clockwise oriented
di-circuits in $D_{\vec M'}$.  Since the symmetric difference of two
perfect matchings consists of disjoint alternating circuits, $D_{\vec
M'}$ arises from $D_M$ by reorienting disjoint di-circuits of $D_M$.
Therefore $\cal S$ is resonant if and only if it is possible to
reorient disjoint di-circuits of $D_{\vec M}$ so that the members of
$\cal S$ will be bounded by clockwise oriented di-circuits.  This is,
in turn, equivalent to requiring that it is possible to reorient
disjoint dicuts of $D\sp *_{\vec M\sp *}$ so that the members of $S$
becomes sink nodes, that is, $S$ is an $F$-stable set of $D\sp *$.
$\bullet $

\medskip By Lemma \ref{equi}, Theorem \ref{Fmax.clar3} when applied to
$D\sp *=(V\sp *,A\sp *)$ and to $F:=\vec M\sp *$ implies the theorem.
$\bullet $ $\bullet $ \medskip

In the last section, we extend the theorem of Abeledo and Atkinson by
deriving a min-max formula, for the maximum number of faces that can
be partitioned into $k$ resonant sets.

\medskip

\section{Optimal $k$-union of sink-stable and $F$-stable sets}


In the preceding section, a result of Gallai was used to prove a
min-max formula for the maximum $w$-weight of the $k$-union of
sink-stable sets and $F$-stable sets for the special case $k=1$.  Now
we solve the case $k\geq 2$ with the help of an extension of Theorem
\ref{Gallaiw}, due to Cameron and Edmonds \cite{Cameron-Edmonds92}.

\THEOREM [Cameron and Edmonds] \label{CE} Let $Q$ be the di-circuits
versus nodes incidence matrix of a digraph $D$.  Let $f:V\rightarrow
{\bf Z}_+\cup \{-\infty \}$ and $g:V\rightarrow {\bf Z}_+\cup \{\infty
\}$ be functions for which $f\leq g$.  The linear system $\{Qx \leq
\widetilde c, \ f\leq x\leq g\}$ is TDI.  $\bullet $ \eT

In the special case when $f\equiv 0$ and $g\equiv 1$, Theorem \ref{CE}
and the l.p. duality theorem immediately gives rise to the following
min-max formula.

\THEOREM \label{CE2} Let $D\sp *=(V,A\sp *)$ be a digraph in which
every node belongs to a di-circuit and let ${\cal K}\sp *$ denote the
set of di-circuits of $D\sp *$.  Let $w:V\rightarrow {\bf Z}_+$ and
$c:A\sp *\rightarrow {\bf Z}_+$ be functions.  Then \eq \max \left
\{\widetilde w(S) :  S\subseteq V, \ \vert S\cap V(K)\vert \leq
\widetilde c(K) \hbox{ for every}\ K\in {\cal K}\sp * \right \}=
\label{(CE2primal)} \eeq \eq \min _{y:{\cal K}\sp * \rightarrow {\bf
Z}_+} \ \left \{\sum _{K\in {\cal K}\sp *} y(K)\widetilde c(K) + \
\sum _{v\in V} \left ( w(v)- \sum _{K\in {\cal K}\sp *, v\in V(K)}
y(K) \right )\sp + \right \}.  \hbox{ $\bullet $ }\ \label{(CE2dual)}
\eeq \eT

We apply this result in the special case when $c$ is the indicator
function of a flat subset of edges.

\THEOREM \label{Fmax} Let $F$ be a flat subset of edges of a strongly
connected digraph $D\sp *=(V,A\sp *)$ with $\vert V\vert \geq 2$ and
let ${\cal K}\sp *$ denote the set of di-circuits of $D\sp *$.  Let
$k\geq 2$ be an integer and $w:V\rightarrow {\bf Z}_+$ an
integer-valued weight-function on the node-set of $D\sp *$.  \ The
maximum $w$-weight of the $k$-union of $F$-stable sets of $D\sp *$ is
equal to \eq \min _{y:{\cal K}\sp *\rightarrow {\bf Z}_+} \ \left \{ k
\sum _{K\in {\cal K}\sp *} y(K) \vert F\cap K\vert + \ \sum _{v\in V}
\left ( w(v)- \sum _{K\in {\cal K}\sp *, v\in V(K)} y(K) \right )\sp +
\right \}.  \label{(Fmaxdual)} \eeq In particular, for a specified
subset $U$ of nodes, the maximum cardinality of the $k$-union of
$F$-stable subsets of $U$ is equal to \eq \min _ { {\cal K}\subseteq
{\cal K}\sp * } \ \left \{ k \sum _{K\in {\cal K}} \vert F\cap K\vert
+ \left \vert U - \cup (V(K):K\in {\cal K})\right \vert \right \}.
\eeq \eT

\Proof Apply Theorem \ref{CE2} to $c':=k\underline \chi _F$ and
observe that \eref{(CE2dual)} transforms to \eref{(Fmaxdual)}.  By
Theorem \ref{Fsinkjel}, a subset $S$ of nodes satisfies the properties
in \eref{(CE2primal)} for $c'$ if and only if $S$ is the $k$-union of
$F$-stable sets.  Hence the first part is a consequence of Theorem
\ref{CE2}.  The second part follows by applying the first one in the
special case $w:=\underline \chi _U$.  $\bullet $

\medskip In the special case when $F$ is a flat transversal of
di-circuits, Theorem \ref{Fmax} is equivalent to a min-max result of
Seb{\H o} \cite{Sebo07} concerning the maximum weight of $k$-unions of
cyclic-stable sets.  Cyclic stability was introduced by Bessy and
Thomasse \cite{Bessy-Thomasse} who proved a min-max result on the
maximum cardinality of a cyclic stable set.  See the last section for
details.

\Corollary [Greene and Kleitman, \cite{Greene-Kleitman}] In a
transitive and acyclic digraph $D'=(U,A')$ the maximum cardinality of
the union of $k$ stable sets is eqaul to $\min \{ k \sum _i \vert
V(P_i)\vert + \vert U- \cup _iV(P_i)\vert :  \{P_1,\dots ,P_q\}$ a set
of disjoint di-paths$\}$.\eCo

\Proof We prove only the non trivial direction $\max\geq \min$.
Extend $D'$ by a new node $z$, add a pair of opposite edges $zu$ and
$uz$ for every $u\in U$.  Let $D=(V,A)$ denote the resulting digraph
and let $F$ be the set of edges entering $z$.  Then $F$ is a flat
transversal of di-circuits of $D$.  Apply the second half of Theorem
\ref{Fmax}.  Observe that each di-circuit of $D$ contains exactly 1
$F$-edge and that the di-circuits in the optimal covering of $U$ can
be made pairwise disjoint in $U$ by the transitivity of $D'$ and hence
their restrictions to $U$ are disjoint di-paths of $D$.  $\bullet $
\medskip

It should be noted that the Greene-Kleitman theorem was derived by
Cameron and Edmonds directly from Theorem \ref{CE}.

We show now how Theorem \ref{Fmax} gives rise to a min-max formula for
the maximum weight of the $k$-union of sink-stable sets.

\THEOREM \label{ksinkmax} Let $D=(V,A)$ be a digraph with no isolated
nodes, $k\geq 2$ an integer, and $w:V\rightarrow {\bf Z}_+$ an
integer-valued weight-function on the node-set of $D$.  Then

$$ \hbox{ $ \max \{ \widetilde w(S):  S\subseteq V$ a $k$-sink-stable
set$\} =$ }\ $$

\eq \min _{ y:  {\cal C}_D\rightarrow {\bf Z}_+} \left \{ k\sum _{C\in
{\cal C}_D} y(C)\eta (C) + \sum _{v\in V} \left (w(v)- \sum _{C\in
{\cal C}_D, v\in V(C)} y(C) \right )\sp + \right \} \label{(ksinkmax)}
\eeq where ${\cal C}_D$ denotes the set of circuits of $D$.  More
concisely, the linear system $\{Qx \leq k\eta , \ 0\leq x\leq \ul 1
\}$ is totally dual integral where $Q$ is the circuit versus node
incidence matrix of $D$ while $\eta $ is a vector the component of
which corresponding to a circuit $C$ is $\eta (C)$.\eT

\Proof Let $S$ be a $k$-sink-stable set.  By Theorem \ref{ksink},
$\vert S\cap V(C)\vert \leq k\eta (C)$ for every circuit $C$ of $D$
and hence $z=\underline \chi _S$ satisfies the primal constraints
$\{Qx\leq k\eta , \ 0\leq x\leq \ul 1\}.$ The trivial direction $\max
\leq \min $ of the l.p. duality theorem implies that $\max \leq \min $
holds in the theorem.

To prove the reverse direction, consider the digraph $D\sp *=(V,A\cup
A')$ arising from $D$ by adding the reverse of each edge of $D$.  Here
$A'$ denotes the set of reverse edges of $D$.  Clearly, $F:=A'$ is a
flat subset since the pair $\{a,a'\}$ is a two-element di-circuit for
every $a\in A$.

\Claim A subset $S$ is a sink-stable set of $D$ if and only if $S$ is
an $F$-stable set of $D\sp *$.  \eCl

\Proof Let $D_2$ be the digraph arising from $D$ by duplicating in
parallel each edge of $D$.  Clearly, $S$ is sink-stable in $D$ if and
only if it is sink-stable in $D_2$.  On the other hand, $D_2$ can be
obtained from $D\sp *$ by reorienting $F=A'$, and hence $S$ is
sink-stable in $D_2$ if and only if it is $F$-stable in $D\sp *$.
$\bullet $

\medskip

It follows from the claim that a subset $S$ of nodes is the $k$-union
of sink-stable sets of $D$ if and only if it is the $k$-union of
$F$-stable sets of $D\sp *.$ By Theorem \ref{Fmax}, the maximum
$w$-weight of a $k$-union of $F$-stable sets of $D\sp *$ is equal to
\eq \min _{y:{\cal K}\rightarrow {\bf Z}_+} \left \{ k \sum _{K\in
{\cal K}\sp *} y(K) \vert F\cap K\vert + \ \sum _{v\in V} \left (
w(v)- \sum _{K\in {\cal K}\sp *, v\in V(K)} y(K) \right )\sp + \right
\} \label{(Fmaxdual2)} \eeq where ${\cal K}\sp *$ denotes the sets of
di-circuits of $D\sp *$

Recall that $D\sp *$ has two types of di-circuits.  A Type I
di-circuits is of form $K_a=\{a,a'\}$ where $a\in A$ while Type II
di-circuits arise from circuits of $D$ by replacing each forward edge
by its reverse or by replacing each backward edge by its reverse.

Consider an optimal (integer-valued) solution $y\sp *$ to
\eref{(Fmaxdual2)} and let $z\sp *(v):= w(v)- \sum [y\sp *(K) :  \
K\in {\cal K}\sp *, v\in V(K) ]$.  Now $z'(v)=z\sp *(v)$ for $v\in
V-\{s,t\}, \ z'(s) = z\sp *(s) + \alpha ,$ and $z'(t) = z\sp *(t) +
\alpha .$ Hence $$\sum _{v\in V} (z'(v))\sp + \leq \sum _{v\in V}
(z\sp *(v))\sp + +2\alpha .$$

Furthermore, $\vert F\cap K_a\vert =1$ implies that $\sum _{K\in {\cal
K}\sp *} y'(K) \vert F\cap K\vert = k \sum _{K\in {\cal K}\sp *} y\sp
*(K) \vert F\cap K\vert -k\alpha $. By combining these observations
with the assumption $k\geq 2$, we obtain that $$k \sum _{K\in {\cal
K}\sp *} y'(K) \vert F\cap K\vert + \ \sum _{v\in V} \left ( w(v)-
\sum _{K\in {\cal K}\sp *, v\in V(K)} y'(K) \right )\sp + \leq $$ $$ k
\sum _{K\in {\cal K}\sp *} y\sp *(K) \vert F\cap K\vert + \ \sum
_{v\in V} \left ( w(v)- \sum _{K\in {\cal K}\sp *, v\in V(K)} y\sp
*(K) \right)\sp +.$$ Since $y\sp *$ is an optimal solution to
\eref{(Fmaxdual2)}, so is $y'$ (and we must have $k=2$), contradicting
the special choice of $y\sp *$.  This contradiction shows that $y\sp
*(K)=0$ for every di-circuit in $D\sp *$ of Type

Suppose now that $K$ is a di-circuit of $D\sp *$ of Type II.  By
reversing the $F$-edges of $K$, we obtain a circuit $C$ of $D$ for
which $\eta (C)\leq \vert F\cap K\vert $. Let $y_0(C):=y\sp *(K)$.  If
$y\sp *(K)>0$ for a di-circuit of $D\sp *$, then we must actually have
$\eta (C) = \vert F\cap K\vert $. Indeed, for if $\eta (C) < \vert
F\cap K\vert $, then $\vert F'\cap K\vert = \eta (C) < \vert F\cap
K\vert $ holds for the reverse di-circuit $K'$ of $K$ and then by
increasing $y\sp *(K)$ with $y\sp *(K)$ and reducing $y\sp *(K)$ to 0
we would obtain a solution to \eref{(Fmaxdual2)} which is better than
the optimal $y\sp *$.  It follows from Theorem \ref{Fmax} that $y_0$
is a solution to \eref{(ksinkmax)} for which $\max \{ w(S) :  S$ a
$k$-sink-set$\} = k\sum _{C\in {\cal C}_D} y_0(C)\eta (C) + \sum
_{v\in V} (w(v)- \sum _{C\in {\cal C}_D, v\in V(C)} y_0(C) \ )\sp +$
from which the min-max result of the theorem follows.  $\bullet $
$\bullet $ \medskip

In the special case $w=\underline \chi _U$, we obtain the following.

\THEOREM Let $D=(V,A)$ be a digraph with no isolated nodes, $k\geq 2$
an integer, and $U\subseteq V$ a prescribed subset of nodes.  The
maximum cardinality of a $k$-sink-stable subset of $U$ is equal to

\eq \min \{k \sum [\eta (C):  C\in {\cal C}] + \vert U-\cup (V(C):
C\in {\cal C})\vert :  {\cal C} \hbox{ a set of circuits}\}.  \hbox{
$\bullet $ }\ \eeq \eT

The same way as the theorem of Abeledo and Atkinson (Theorem
\ref{Abeledo-Atkinson}) was derived from Theorem \ref{Fmax.clar3}, the
following result can be obtained from the second part of Theorem
\ref{Fmax}.

\THEOREM Let $G=(S,T;E)$ be a 2-connected perfectly matchable
bipartite plane graph and $k\geq 2$ an integer.  The maximum number of
faces that can be partitioned into $k$ resonant sets is equal to the
minimum of $$ \hbox{ $k \sum _i$ val$(B_i) + $the number of faces
avoided by each $B_i$ }\ $$ where the minimum is taken over all
choices of feasible cuts $B_1,\dots ,B_q$.  In particular, the set of
faces of $G$ can be partitioned into $k$ resonant sets if and only if
for any set ${\cal B}$ of feasible cuts the number of faces
intersected by ${\cal B}$ is at most $k$ times the total value of
${\cal B}$.  $\bullet $ \eT

\Corollary The faces of a 2-connected perfectly matchable bipartite
plane graph can be partitioned into $k$ resonant sets if and only if,
for every feasible cut $B$ of $G$, the number of faces intersected by
$B$ is at most $k${\em val}$(B)$.  $\bullet $ \eCo

\medskip

\section{Link to cyclic stable sets}

Suppose that $D=(V,A)$ is a strongly connected loopless digraph on
$n\geq 2$ nodes and consider a linear order ${\cal L}=[v_1,\dots
,v_n]$ of the nodes of $D$.  An edge $e$ of $D$ is a forward edge if
its tail precedes its head, otherwise $e$ is a backward edge.

Let $P$ be a regular $n$-gon in a horizontal plane and assign the
nodes of $V$ to the vertices of $P$ in this order.  In this way, we
arrive at a cyclic order ${\cal O}=(v_1,\dots ,v_n)$ of $D$.  A set of
consecutive elements is called an {\bf interval} of $\cal O$.  For
example, both $\{v_2,v_3,v_4\}$ and $\{v_{n-1},v_n,v_1,v_2\}$ are
intervals.  Each edge $uv$ of $D$ can be represented in the plane by
an arc going clockwise outside $P$.  Clearly, the linear order
$[v_i,\dots ,v_n, v_1,\dots ,v_{i-1}]$ defines the same cyclic order
for each $v_i$.  Each of these $n$ linear orders is called an opening
of ${\cal O}$.  The backward edges of the opening ${\cal
L}'=[v_i,\dots ,v_n, v_1,\dots ,v_{i-1}]$ of ${\cal O}$ is called the
edge-set {\bf belonging} to the opening.

Let $K$ be a di-circuit of $D$.  Starting from a node $v$ of $K$ and
going along $K$ we arrive back to $v$.  In the plane, this simple
closed walk goes around $P$ one or more times.  This number is called
the winding number or the index of $K$ and is denoted by ind$(K)$.  It
follows from this definition that if $F$ denotes the set of edges
belonging to an opening of $\cal O$, then \eq \hbox{ ind}(K) = \vert
F\cap K\vert . \label{(indF)} \eeq

For example, if $D$ itself is a di-circuit $K$ consisting of the edges
$\{v_1v_2, \dots ,v_{n-1}v_n, v_nv_1\}$, then the index of $K$ with
respect to the cyclic order $(v_1,\dots ,v_n)$ is 1 while ind$(K)=n-1$
with respect to the reverse cyclic order $(v_n,\dots ,v_1).$

These notions were introduced by Bessy and Thomass\'e
\cite{Bessy-Thomasse} who called a cyclic order of $D$ {\bf coherent}
if each edge of $D$ belongs to a di-circuit of index 1. They proved
that every strong digraph has a coherent ordering.  Let $\cal O$ be a
cyclic ordering and $F$ the set of edges belonging to an opening of
$\cal O$.

Iwata and Matsuda \cite{Iwata-Matsuda} observed the following link
between flat transversals of di-circuits and coherent cyclic
orderings.

\Lemma [Iwata and Matsuda] \label{flat-coherent} Let $D=(V,A)$ be a
strongly connected digraph.  A subset $F$ of edges is a flat
transversal of di-circuits if and only of $F$ belongs to an opening of
a coherent ordering of $D$.  \eL

\Proof (Outline) If $F$ belongs to an opening of a cyclic order ${\cal
O}$, then $F$ is clearly a transversal of di-circuits.  If ${\cal O}$
is, in addition, coherent, that is, if each edge belongs to a
di-circuit of index 1, then $F$ is flat since ind$(K)= \vert F\cap
K\vert $ for every di-circuit.

Conversely, if $F$ is a flat transversal of di-circuits, then $F$ is
certainly a minimal transversal with respect to inclusion.  An easy
excercise shows that the digraph $D_F$ arising from $D$ by reversing
the elements of $F$ is acyclic.  Hence any topological ordering $\cal
L$ of $D_F$ has the property that the elements of $F$ (in $D$) are
precisely the backward edges.  Therefore the cyclic order determined
by $\cal L$ is coherent.  $\bullet $

\medskip Due to this correspondence, the existence of a coherent
cyclic order is equivalent to Knuth's lemma on the existence of a flat
transversal of di-circuits.

Bessy and Thomasse called the exchange of two consecutive elements $u$
and $v$ in a cyclic order {\bf elementary} if there is no edge (in
either direction) between $u$ and $v$.  They called two cyclic orders
equivalent if one can be obtained from the other by a sequence of
elementary exchanges.  Finally, a stable set of nodes is {\bf cyclic
stable} with respect to a given cyclic order $\cal O$ if there is an
equivalent cyclic order where $S$ forms an interval.

From a complexity point of view, a slight disadvantage of this
definition is that it does not show (as it is) that cyclic stability
is an NP-property.  Indeed, in principle it could be the case that a
cyclic order can be obtained from an equivalent cyclic order only by a
sequence of exponentially many elementary exchanges, and in such a
case the definition would not provide a polynomally checkable
certificate for cyclic stability.  A. Seb{\H o} \cite{Sebo08},
however, pointed out that a cyclic ordering can always be obtained
from an equivalent cyclic order by a sequence of at most $n\sp 2$
elementary exchanges.  Hence cyclic stability is an NP-property.
Moreover, Seb{\H o} proved the following co-NP characterization of
cyclic stability.

\THEOREM [\cite{Sebo07}, Statement (5)] \label{cyclejel} A subset $S$
of nodes of a strongly connected digraph $D$ is cyclic stable with
respect to a coherent cyclic ordering if and only if $\vert S\cap
V(K)\vert \leq $ind$(K)$ for every di-circuit $K$ of $D$.  $\bullet $
\eT

The proof of this theorem provides a polynomial algorithm that either
finds a sequence of elementary exchanges that transform $S$ into an
interval or else it finds a di-circuit $K$ violating the inequality in
the theorem.

The two main theorems of Bessy and Thomass\'e \cite{Bessy-Thomasse}
are as follows.

\THEOREM [Bessy and Thomass\'e] \label{BT1} Given a strong digraph
$D=(V,A)$ along with a coherent cyclic ordering, the maximum
cardinality of a cyclic stable set of $D$ is equal to the minimum
total index of di-circuits covering $V$.  $\bullet $ \eT

\THEOREM [Bessy and Thomass\'e] \label{BT2} Let $D=(V,A)$ be a strong
digraph along with a coherent cyclic ordering and let $k\geq 2$ be an
integer.  The node-set of $D$ can be partitioned into $k$ cyclic
stable sets if and only if $\vert K\vert \leq k\ ${\em ind}$(K)$ for
every di-circuit $K$.  $\bullet $ \eT

As a common generalization of the theorems of Bessy and Thomass\'e,
Seb{\H o} \cite{Sebo07} proved the following.

\THEOREM [Seb{\H o}] \label{Sebo} Let $D=(V,A)$ be a strong digraph
along with a coherent cyclic ordering.  Let $k\geq 1$ be an integer
and $U$ a subset of nodes.  The maximum cardinality of the union of
$k$ cyclic stable sets of $D$ is equal to $\min \{k \sum _i$ {\em
ind}$(K_i) + \vert U -\cup _iV(K_i)\vert :  \{K_1,\dots ,K_q\}$ a set
of di-circuits$\}$.  \eT

Seb{\H o} actually proved this result in a more general form by
providing a min-max formula for the maximum $w$-weight of the
$k$-union of cyclic stable sets.

By combining Theorems \ref{Fsinkjel} and \ref{cyclejel}, we obtain by
\eref{(indF)} the following.

\Lemma \label{link} Let $F$ be a flat transversal of di-circuit of a
strong digraph $D$ and let ${\cal O}=(v_1,\dots ,v_n)$ be a coherent
cyclic order so that $F$ belongs to an opening of $\cal O$.  Then a
subset $S$ of nodes is $F$-stable if and only if $S$ is cyclic stable.
$\bullet $ \eL

This lemma implies that Theorem \ref{Sebo} is equivalent to the that
special case of the second half of Theorem \ref{Fmax} when $F$ is not
only flat but it is a transversal of di-circuits, as-well.  Note that
requiring only the flatness of $F$ in Theorem \ref{Fmax} allowed us to
derive the theorem of Abeledo and Atkinson and its extension.\\
\bigskip

{\bf Acknowledgements} \ Many thanks are due to Zolt\'an Kir\'aly for
the valuable discussions in an early phase of this research.  The
second author received a grant (no.  CK 80124) from the National
Development Agency of Hungary, based on a source from the Research and
Technology Innovation Fund.  Part of research was done while he
visited the Research Institute for Mathematical Sciences, Kyoto
University, 2008, and the Institute of Discrete Mathematics,
University of Bonn, 2011.

\pagebreak

 \end{document}